\documentclass[11pt]{article}
\usepackage{amssymb, amsthm, amsmath, amscd}
\newtheorem{thm}{Theorem}[section]
\newtheorem{lem}[thm]{Lemma}

\newtheorem{cor}[thm]{Corollary}
\newtheorem{NN}[thm]{}
\theoremstyle{definition}\newtheorem{df}[thm]{Definition}
\theoremstyle{definition}
\theoremstyle{definition}

\renewcommand{\phi}{\varphi}

\newcommand{\Z}{\mathbb{Z}}
\newcommand{\Q}{\mathbb{Q}}

\newcommand{\C}{\mathbb{C}}
\newcommand{\T}{\mathbb{T}}

\newcommand{\hm}{homomorphism}

\newcommand{\andeqn}{\,\,\,{\rm and}\,\,\,}
\newcommand{\rforal}{\,\,\,{\rm for\,\,\,all}\,\,\,}
\newcommand{\CA}{$C^*$-algebra}

\newcommand{\af}{{\alpha}}
\newcommand{\bt}{{\beta}}

\newcommand{\beq}{\begin{eqnarray}}
\newcommand{\eneq}{\end{eqnarray}}
\newcommand{\tforal}{\,\,\,\text{for\,\,\,all}\,\,\,}
\newcommand{\tand}{\,\,\,\text{and}\,\,\,}
\newcommand{\p}{\mathfrak{p}}
\newcommand{\q}{\mathfrak{q}}

\title{Localizing the Elliott Conjecture at Strongly  Self-absorbing \CA s, II\\
-----An Appendix }
\author{Huaxin Lin\\
Department of Mathematics\\
East China Normal University,\\
 Shanghai, China\\
 and \\
 Department of Mathematics\\
 University of Oregon\\
 Eugene, Oregon, 97403\\
 U.S.A
 }
\date{}

\begin{document}

\maketitle

\begin{abstract}
This note provides some technical support to the proof of a result
of W. Winter which shows that two unital separable simple amenable
${\cal Z}$-absorbing \CA s with locally finite decomposition
property satisfying the UCT whose projections separate the traces
are isomorphic if their $K$-theory is  finitely generated and
their Elliott invariants are the same.

\end{abstract}

\section{Introduction}
In \cite{W}, W. Winter provided a fascinating method for the
Elliott program of classification of amenable \CA s. Let $A$ and
$B$ be two unital separable simple amenable \CA s which are ${\cal
Z}$-absorbing. Winter showed that if $A\otimes C$ is isomorphic to
$B\otimes C,$ for any unital UHF-algebra $C,$  then there is a way to
show that $A$ is actually isomorphic to $B.$ It is known that for
many separable amenable simple \CA s, $A\otimes C$  has many known
properties that $A$ does not have. For example, it was proved
(also by Winter) that, for every ${\cal Z}$-absorbing \CA\, $A$
with locally finite decomposition rank whose projections separate
the traces, $A\otimes C$ has tracial rank zero for every
UHF-algebra $C$ (see \cite{W1}). However, $A$ itself may not have
finite tracial rank. Winter's method shows that two unital
separable simple amenable ${\cal Z}$-absorbing \CA s with locally
finite decomposition rank and with finitely generated $K$-theory
whose projections separate the traces are isomorphic if their
Elliott invariants are isomorphic. The main purpose of this note is
to provide some technical support to the proof of the above
mentioned result of Winter.

\vspace{0.2in}

{\bf Acknowledgments}  This work was done when the author was in
East China Normal University during the summer of 2007. He
acknowledges the support from Shanghai Priority Academic
Disciplines.

\section{Some notation}

\begin{NN}\label{Mtorus}
{\rm Let $A$ and $B$ be two unital \CA s. Suppose that $\phi, \psi:
A\to B$ are two \hm s. Define the mapping torus of $\phi$ and $\psi$
as follows:
\beq\label{Dmt-1}
M_{\phi, \psi}=\{x\in C([0,1], B): x(0)=\phi(a)\andeqn x(1)=\psi(a)
\,\,\,\text{for some}\,\,\, a\in A\}.
\eneq
  Thus one
obtains an exact sequence:
\beq\label{Dmt-2}
0\to SB \stackrel{\imath}{\to} M_{\phi, \psi}\stackrel{\pi_0}{\to}
A\to 0,
\eneq
where $\pi_0: M_{\phi,\psi}\to A$ is the point-evaluation at $t=0.$

 Suppose that $A$ is a separable amenable \CA.  From
(\ref{Dmt-2}), one obtains an element in $Ext(A,SB).$ In this case
we identify $Ext(A,SB)$ with $KK^1(A,SB)$ and $KK(A,B).$

Suppose that $[\phi]=[\psi]$ in $KL(A,B).$ The mapping torus
$M_{\phi, \psi}$ corresponds to a trivial element in $KL(A, B).$ It
follows that there are two exact sequences:
\beq\label{Dmt-3}
&&0\to K_1(B)\stackrel{\imath_*}{\to} K_0(M_{\phi,
\psi})\stackrel{(\pi_0)_*}{\to} K_0(A)\to 0\tand\\\label{Dmt-4}
&&0\to K_0(B)\stackrel{\imath_*}{\to}K_1(M_{\phi, \psi})
\stackrel{(\pi_0)_*}{\to} K_1(A)\to 0.
\eneq
which are pure extensions of abelian groups.

 }

\end{NN}

\begin{df}\label{Dr}
{\rm  Now let $B$ be a unital \CA\, with non-empty tracial state
space $T(B).$ Let $u\in M_l(M_{\phi, \psi})$ be a unitary which is
a piecewise smooth function on $[0,1].$ For each $\tau\in T(B),$
we denote by $\tau$ the trace $\tau\otimes Tr$ on $M_l(B),$
where $Tr$ is the standard trace on $M_l.$  Define
\beq\label{Dr-1}
R_{\phi,\psi}(u)(\tau)={1\over{2\pi i}}\int_0^1
\tau({du(t)\over{dt}}u(t)^*)dt.
\eneq

When $[\phi]=[\psi]$ in $KL(A,B)$ and $\tau\circ \phi=\tau\circ
\psi$ for all $\tau\in T(B),$ there exists a \hm\,
$$
R_{\phi, \psi}: K_1(M_{\phi, \psi})\to Aff(T(B))
$$
defined by
$$
R_{\phi, \psi}([u])(\tau)={1\over{2\pi i}}\int_0^1
\tau({du(t)\over{dt}}u(t)^*)dt.
$$

 We will call $R_{\phi, \psi}$ the {\it rotation map} for the pair
$\phi$ and $\psi.$

}
\end{df}

Moreover,  the following diagram commutes:
$$
\begin{array}{ccccc}
K_0(B) && \stackrel{\imath_*}{\longrightarrow} && K_1(M_{\phi, \psi})\\
& \rho_B\searrow && \swarrow R_{\phi,\psi} \\
& & Aff(T(B)) \\
\end{array}
$$

See section 3 of \cite{Laut} for more information.

\vspace{0.2in}

\begin{df}
If furthermore, $[\phi]=[\psi]$ in $KK(A,B)$ and $A$ satisfies the
Universal Coefficient Theorem, using Dadarlat-Loring's notation, one
has the following splitting exact sequence:
\beq\label{eta-1-}
0\to \underline{K}(SB)\,{\stackrel{[\imath]}{\to}}\,
\underline{K}(M_{\phi_1,\phi_2})\,{\stackrel{[\pi_0]}{\rightleftarrows}}_{\theta}
\,\,\underline{K}(A)\to 0.
\eneq
In other words there is  $\theta\in Hom_{\Lambda}(\underline{K}(A),
\underline{K}(M_{\phi_1,\phi_2}))$ such that $[\pi_0]\circ
\theta=[\rm id_A].$ In particular, one has a monomorphism
$\theta|_{K_1(A)}: K_1(A)\to K_1(M_{\phi, \psi})$ such that
$[\pi_0]\circ \theta|_{K_1(A)}=({\rm id}_A)_{*1}.$ Thus, one may
write
\beq\label{eta-1}
K_1(M_{\phi, \psi})=K_0(B)\oplus K_1(A).
\eneq
Suppose also that $\tau\circ \phi_1=\tau\circ \phi_2$ for all
$\tau\in T(B).$  Then  one obtains the \hm\,
\beq\label{eta-2}
R_{\phi,\psi}\circ \theta|_{K_1(A)}: K_1(A)\to Aff(T(B)).
\eneq

We write $${\tilde \eta}_{\phi_1,\phi_2}=0,$$ if $R_{\phi,\psi}\circ
\theta=0,$ i.e., $\theta(K_1(A))\in {\rm ker}R_{\phi_1,\phi_2}$ for
some such $\theta.$ Thus, $\theta$ also gives the following:
$$
{\rm ker}R_{\phi,\psi}={\rm ker}\rho_B\oplus K_1(A).
$$

\end{df}

\begin{df}
Let $A$ and $B$ be two unital \CA s and let $\phi, \psi: A\to B$
be two \hm s. We say $\phi$ and $\psi$ are {\it asymptotically
unitarily equivalent} if there exists a continuous path of
unitaries $\{u(t):t\in [0, 1)\}$ of $B$ such that
\beq\label{d-1}
\lim_{t\to 1}{\rm ad}\, u(t)\circ \phi(a)=\psi(a)\rforal a\in A.
\eneq
We say $\phi$ and $\psi$ are {\it strongly asymptotically
unitarily equivalent} if $u(t)$ can be so chosen that $u(0)=1.$

\end{df}

We use the following result in the proof.

\begin{thm}{\rm (Theorem 9.1 of \cite{Laut})}\label{T0}
Let $A$ be a unital AH-algebra and let $B$ be a unital simple
\CA\, with tracial rank zero. Suppose that $\phi_1, \phi_2: A\to
B$ are two monomorphisms. Then $\phi_1$ and $\phi_2$ are
asymptotically unitarily equivalent if and only if
\beq
[\phi_1]=[\phi_2]\,\,\,{\rm in}\,\,\, KK(A,B),\,\,
\tau\circ\phi_1=\tau\circ \phi_2\tforal \tau\in
T(A)\,\,\,and\,\,\, {\tilde \eta}_{\phi_1,\phi_2}=\{0\}.
\eneq
\end{thm}

\vspace{0.2in}

 In what follows, $Q$ denotes the UHF-algebra with
$K_0(Q)=\Q$ and $[1_Q]=1,$ and if $\p$ is a supernatural number
$M_\p$ denotes the UHF-algebra associated with the supernatural
number $\p.$

\section{The Main results}

\begin{lem}\label{L00}
Let $A$ be a unital \CA\, and let $h_1,h_2,...,h_n$ be
self-adjoint elements in $A.$ Suppose that $v$ is any unitary in
$A$ and
$$
u(t)=v\prod_{j=1}^n e^{ih_jt}\, \,\,\,\,\,t\in [0,1].
$$
Then,
$$
\int_0^1\tau({du(t)\over{dt}}u(t)^*)dt=\sum_{j=1}^n\tau(h_j)
$$
for all $\tau\in T(A).$
\end{lem}

\begin{proof}
Note that for any unitary $w\in A$ and any tracial state $\tau\in
T(A),$
$$
\tau(wh_jw^*)=\tau(h_j),\,\,\,j=1,2,...,n.
$$
It follows that
$$
\tau({du(t)\over{dt}}u(t)^*)= \sum_{j=1}^n\tau(h_j)
$$
for all $\tau\in T(A).$ Thus the lemma follows.

\end{proof}

\begin{lem}\label{l1}
Let $B$ be a unital separable simple amenable \CA\, such that
$B\otimes Q$ has tracial rank zero. Let $A$ be a unital separable
amenable simple \CA\, with tracial rank zero satisfying the UCT
such that $K_i(A)$ is torsion free ($i=0,1$).

Suppose that $\phi_1, \phi_2: A\to B\otimes Q$ are two unital \hm
s with
$$
[\phi_1]=[\phi_2]\,\,\,{\rm in}\,\,\, KK(A, M\otimes Q).
$$
Suppose that $\phi_1$ induces an affine homeomorphism
$(\phi_1)_{\sharp}: T(B\otimes Q)\to T(A)$ by
$$
(\phi_1)_{\sharp}(\tau)(a)=\tau\circ \phi_1(a)
$$
for all $\tau\in T(B\otimes Q)$ and $a\in A.$

Then there exists an automorphism $\af\in Aut(\phi_1(A))$ with
$[\af]=[{\rm id}_{\phi_1(A)}]$ in $KK(\phi_1(A),\phi_1(A))$ such
that $ \af\circ \phi_1 $ and $\phi_2$ are strongly asymptotically
unitarily equivalent.
\end{lem}

\begin{proof}
Let
$$
M_{\phi_1,\phi_2}=\{f\in C([0,1], B\otimes Q): f(0)=\phi_1(a)
\andeqn f(1)=\phi_2(a)\,\,\,{\rm for\,\,\,some}\,\,\,a\in A\}.
$$
Note that, since $K_i(B\otimes Q)$ is torsion free, $i=0,1,$
$$
KK(A, B\otimes Q)\cong Hom(K_*(A), K_*(B\otimes Q)).
$$
 Since we assume that $[\phi_1]=[\phi_2]$ in
$KK(A, B\otimes Q),$ there exists a \hm\, $\theta: K_i(A)\to
K_i(M_{\phi_1,\phi_2})$ such that
$$
0\to K_{i-1}(B\otimes Q)\to K_i(M_{\phi_1,
\phi_2}){\stackrel{[\pi_0]}{\rightleftarrows}}_{\theta} \,\,\,
K_i(A)\to 0
$$
splits ($i=0,1$).

Since $A$ has real rank zero, we also have
\beq\label{tr1}
\tau\circ \phi_1(a)=\tau\circ \phi_2(a)\tforal a\in A
\eneq
and for all $\tau\in T(B\otimes Q).$

Let $R_{\phi_1,\phi_2}: K_1(M_{\phi_1,\phi_2})\to Aff(T(B\otimes
Q))$ be the rotation map. Put $C=\phi_1(A).$
By the classification theorem of \cite{Lnduke}, $C$ is a unital simple AH-algebra with no dimension
growth and with real rank zero. Since $K_i(C)$ is torsion free, it follows that $C$ is a unital simple
$A\T$ -algebra with real rank zero.

By 4.1 (see also Theorem 4.4) of \cite{KK1}, there exists an
automorphism $\af\in Aut(C)$ with $[\af]=[{\rm id}_{C}]$ in
$KK(C,C)$ satisfying the following: there is $\theta_1: K_1(C)\to
K_1(M_{\af})$ so that
$$
(R_{\af}\circ
\theta_1)((\phi_1)_{\sharp}(\tau))=-(R_{\phi_1,\phi_2}\circ
\theta)(\tau)\tforal \tau\in T(B\otimes Q),
$$
where
$$
M_\af=\{g\in C([0,1], C): g(0)=a\andeqn g(1)=\af(a)\,\,\,{\rm
for\,\,\,some }\,\,\,a\in C\}
$$
and
$$
R_\af(u)(\tau)={1\over{2\pi}}\int_0^1
\tau({du(t)\over{dt}}u(t)^*)dt
$$
for all unitaries $u\in M_{k}(M_\af)$ ($k=1,2,...,$) and $\tau\in
T(C).$ Note that $[\af]=[{\rm id}_C]$ in $KK(C,C).$ Therefore, one
computes that
$$
[\af\circ\phi_1]=[\phi_1]\,\,\,{\rm in}\,\,\,KK(A, B\otimes Q).
$$

 Now consider the mapping torus
$$
M_{\af\circ \phi_1,\phi_2}=\{f\in C([0,1], B\otimes Q):
f(0)=\af\circ \phi_1(a)\andeqn f(1)=\phi_2(a)\,\,\,{\rm
for\,\,\,some}\,\,\,a\in A\}.
$$

 Define $\theta_2: K_1(C)\to K_1(M_{\af\circ \phi_1,
 \phi_2})$ as follows:

 Let $k>0$ be an integer and $u\in M_k(M_{\af\circ \phi_1,
 \phi_2})$ be a unitary.

 We may assume that there is a unitary $w(t)\in M_k(M_{\phi_1, \phi_2})$ such that
 \beq\label{path1}
 w(0)=\phi_1(u'),w(1)=\phi_2(u'),[u']=[u]\,\,\,{\rm in}\,\,\, K_1(A)\\
 \andeqn
 \theta([u])=[w(t)]\,\,\,{\rm in}\,\,\,K_1(M_{\phi_1,\phi_2})
 \eneq
 for some unitary $u'\in M_k(A).$
To simplify notation, without loss of generality, we may assume
that there are $h_1,h_2,...,h_n\in M_k(A)_{s.a.}$ such that
$$
u^*u'=\prod_{j=1}^n\exp(ih_j).
$$
Define $z(t)=u\prod_{j=1}^n\exp(h_jt)$ ($t\in [0,1]$). Consider
$\{\phi_1(z(t)): t\in [0,1]\}.$ Then
$$
\phi_1(z(0))=\phi_1(u)\andeqn \phi_1(z(1))=\phi_1(u').
$$
Moreover, by \ref{L00},
$$
\int_0^1\tau({d\phi_1(z(t))\over{dt}}\phi_1(u(t))^*)dt=\sum_{j=1}^n\tau(\phi_1(h_j))
$$
for all $\tau\in T(A).$ Consider $Z(t)=\phi_2(z(1-t)).$ Then
$$
Z(0)=\phi_2(u')\andeqn Z(1)=\phi_2(u).
$$
$$
\int_0^1\tau({d\phi_2(z(1-t))\over{dt}}\phi_2(z(1-t))^*)dt=-\sum_{j=1}^n\tau(\phi_2(h_j))
$$
for all $\tau\in T(A).$ Note that
$$
\tau(\phi_2(h_j))=\tau(\phi_1(h_j)),\,\,\,j=1,2,...,n.
$$
It follows that
$$
\int_0^1\tau({d\phi_1(z(t))\over{dt}}\phi_1(u(t))^*)dt
+\int_0^1\tau({d\phi_2(z(1-t))\over{dt}}\phi_2(z(1-t))^*)dt=0
$$
for all $\tau\in T(A).$

Therefore,
 without loss of generality, we may assume that $u=u'$ in
(\ref{path1}). We may also assume that both paths are piecewise
smooth.

We may also assume that there is a unitary $s(t)\in
  M_k(M_{\af\circ \phi_1, \phi_1})$ such that
  \beq
 s(0)=\af\circ\phi_1(u),s(1)=\phi_1(u)\,\,\,{\rm in}\,\,\, K_1(A)\\
 \andeqn
 \theta_1([u])=[s(t)]\,\,\,{\rm in}\,\,\,K_1(M_{\af\circ\phi_1,\phi_1}).
 \eneq

 Define
 $\theta_2([u])=[v],$ where
 \beq\label{path}
v(t)=\begin{cases} s(2t)\,\,\,\text{if}\,\,\, t\in
[0,1/2)\\
w(2(t-1/2))\,\,\,\text{if}\,\,\, t\in [1/2,1],
\end{cases}
\eneq
Thus $\theta_2$ gives a \hm\, from $K_1(A)$ to $K_1(M_{\af\circ
\phi_1,\phi_2})$ such that $(\pi_0)_{*1}\circ \theta_2={\rm
id}_{K_1(A)}.$ Since $[\af\circ \phi_1]=[\phi_1]$ in $KK(A,C)),$
we also have
$$
[\af\circ \phi_1]=[\phi_2]\,\,\,{\rm in}\,\,\,KK(A, B\otimes Q).
$$
There is also a \hm\, $\theta_2': K_0(A)\to K_0(M_{\af\circ
\phi_1,\phi_2})$ such that $(\pi_0)_{*0}\circ \theta_2'=[{\rm
id}_{K_0(A)}].$

Then
\beq
R_{\af\phi_1, \phi_2}([u])(\tau)&=& {1\over{2\pi}}\int_0^1
\tau({dv(t)\over{dt}}v(t)^*)dt\\
&=&{1\over{2\pi}}\int_0^{1/2} \tau({ds(2t)\over{dt}}s(2t)^*)dt+\\
&&{1\over{2\pi}}\int_{1/2}^{1}
\tau({dw(2(t-1/2))\over{dt}}w(2(t-1/2))^*)dt\\
&=& R_{\af\circ \phi_1,\phi_1}\circ
\theta_1([u])(\phi_{\sharp}(\tau))+R_{\phi_1, \phi_2}\circ
\theta([u])(\tau)=0
\eneq
for all $\tau\in T(B\otimes Q).$

Thus ${\tilde{\eta}_{\af\circ \phi_1, \phi_2}}=0.$ Note that, by
\cite{Lnduke}, $A$ is an AH-algebra. It follows from Theorem
\ref{T0} that $\af\circ \phi_1$ and $\phi_2$ are asymptotically
unitarily equivalent. Since $K_1(B\otimes Q)$ is
divisible,\linebreak
 $H_1(K_0(B\otimes Q), K_1(B\otimes
Q))=K_1(B\otimes Q),$ and by 11.5 of \cite{Laut}, we conclude that
$\af\circ \phi_1$ and $\phi_2$ are strongly asymptotically
unitarily equivalent.

\end{proof}

\begin{lem}\label{2L1}
Let $B$ be a unital separable simple amenable \CA\, such that
$B\otimes Q$ has tracial rank zero. Let $\p$ be a supernatural
number of infinite type and let $A=C\otimes M_{\mathfrak{p}}$ be a
unital separable amenable simple \CA\, with tracial rank zero
satisfying the UCT such that $K_i(A)=Tor(K_i(A))\oplus G_i,$ where
$G_i$ is torsion free ($i=0,1$). Suppose that $\phi_1, \phi_2:
A\to B\otimes Q$ are two unital \hm s with
$$
[\phi_1]=[\phi_2]\,\,\, in\,\,\, KK(A, M\otimes Q).
$$
Suppose also that $\phi_1$ induces an affine homeomorphism
$(\phi_1)_{\sharp}: T(B\otimes Q)\to T(A)$ by
$$
(\phi_1)_{\sharp}(\tau)(a)=\tau\circ \phi_1(a)
$$
for all $\tau\in T(B\otimes Q)$ and $a\in A.$

Then there exists an automorphism $\af\in Aut(\phi_1(A))$ with
$[\af]=[{\rm id}_{\phi_1(A)}]$ in $KL(\phi_1(A),\phi_1(A))$ such
that $ \af\circ \phi_1 $ and $\phi_2$ are strongly asymptotically
unitarily equivalent.
\end{lem}

\begin{proof}
Let $d: K_0(A)\to Aff(T(A))$ be the \hm\, induced by $\tau([p])$
for all $\tau\in T(A)$ and projections $p\in A.$ There is a \hm\,
$h_1: M_{\mathfrak{p}}\to M_{\mathfrak{p}}$ such that $h_1(1)=e$
for some projection $e\in M_{\mathfrak{p}}$ with $e\not=0$ and
$e\not=1.$ This $h_1$ gives an injective \hm\, $\gamma_1:
d(K_0(A))\to d(K_0(A))$ such that $\gamma_1(r)=[e]r$ for $r\in
d(K_0(C\otimes M_{\mathfrak{p}})).$ Put $A_1=\phi_1(A).$
 Define $h\in Hom(K_*(A_1), K_*(A_1))$ such that
$$
h|_{K_0(A_1)}=\gamma_1\circ d, , h|_{G_1}={\rm id}|_{G_{1}}\andeqn
$$
$$
h|_{Tor(K_1(A_1))}=\{0\}.
$$
There is a unital simple $A\T$-algebra $D$ with real rank zero
such that
$$
(K_0(D), K_0(D)_+, [1_D], K_1(D))=(h(K_0(A_1)), h((K_0(A_1))_+),
h([1_{A_1}]), G_1).
$$
So $K_i(D)=G_i$ ($i=0,1$). It follows from \cite{Lnduke} that
there exists a unital \hm\, $\psi_1': A_1\to D$ such that
$(\psi_1')_{*i}=h.$ Moreover, there is a \hm\, $\imath: D\to A_1$
so that $(\imath)_{*i}={\rm id}_{G_i},$ $i=0,1.$

Put
$$
\kappa=[{\rm id}_{A_1}]-[\imath\circ \psi_1']\in KL(A_1, A_1).
$$
Note that $\kappa\in KL(A_1, A_1)_{+}$ (see \cite{Lnduke}). It
follows from \cite{Lnduke} that there is a \hm\, $\psi_2: A_1\to
A_1$ such that $\psi_2(1_{A_1})=1_{A_1} -\imath\circ
\psi_1'(1_{A_1})=e_1$ (for some projection $e_1\in A_1$) and
$[\psi_2]=\kappa.$  Define
$$
\Phi(a)=\imath\circ \psi_1'(a)+\psi_2(a)\tforal a\in A_1.
$$

It is  clear that
$$
[\Phi]=[{\rm id}_{A_1}]\,\,\,{\rm in}\,\,\,KL(A_1, A_1).
$$
It follows that
$$
[\Phi\circ \phi_2]=[\phi_1]\,\,\,{\rm in}\,\,\,KK(A, B\otimes Q).
$$

 Let $\theta_1: K_i(A)\to K_i(M_{\Phi\circ\phi_1, \phi_2})$
($i=0,1$) be such that
$$
0\to K_{i-1}(B\otimes Q)\to K_i(M_{\Phi\circ\phi_1,
\phi_2}){\stackrel{[\pi_0]}{\rightleftarrows}}_{\theta_1} \,\,\,
K_i(A)\to 0
$$
splits, where
$$
M_{\Phi\circ \phi_1, \phi_2}=\{ f\in C([0,1], B\otimes Q):
f(0)=\Phi\circ \phi_1(a)\andeqn f(1)=\phi_2(a)\,\,\, {\rm
for\,\,\,some}\,\,\, a\in C\}.
$$
Let $\eta_1=R_{\Phi\circ \phi_1,\phi_2}\circ \theta_1.$ We will
identify $D$ with $\imath(D)$ and we  identify $Aff(T(D))$ with
$Aff(T(A)).$ It follows from \cite{KK1} that there exists an
automorphism $\af_1\in Aut(D)$ with $[\af_1]=[{\rm id}_D]$ in
$KK(D,D)$ satisfying the following: there is $\theta_2: K_1(D)\to
K_1(M_{\af_1})$ such that
$$
R_{\af_1}\circ
\theta_2((\phi_1)_{\sharp}(\tau))=-\eta_1(\tau)\tforal \tau\in
T(B\otimes Q),
$$
where
$$
M_{\af_1}=\{f\in C([0,1], D):f(0)=\af_1(c)\andeqn f(1)=c\,\,\,{\rm
for\,\,\,some}\,\,\, c\in D\}.
$$

Now define $\af': D\oplus e_1\Phi(A_1)e_1\to D\oplus
e_1\Phi(A_1)e_1$  by
$$
\af'(c)=\af_1(c)\tforal c\in D\andeqn \af'(a)=a\tforal a\in
e_1\Phi(A_1)e_1.
$$
Define $\af=\af'\circ \Phi: A_1\to A_1.$ Then
$$
[\af]=[{\rm id}_{A_1}]\,\,\,{\rm in}\,\,\, KL(A_1,A_1)
$$
and
$$
[\af\circ \phi_1]=[\phi_2]\,\,\,{\rm in}\,\,\, KK(A, B\otimes Q).
$$
 Let $k>0$ be an integer and $u\in M_k(M_{\Phi\circ \phi_1,
 \phi_2})$ be a unitary.

 We may assume that there is a unitary $w(t)\in M_k(M_{\phi_1, \phi_2})$ such that
 \beq\label{path1+1}
 w(0)=\Phi\circ\phi_1(u'),w(1)=\phi_2(u'),[u']=[u]\,\,\,{\rm in}\,\,\, K_1(A)\\
 \andeqn
 \theta_1([u])=[w(t)]\,\,\,{\rm in}\,\,\,K_1(M_{\Phi\circ\phi_1,\phi_2})
 \eneq
 for some unitary $u'\in M_k(A).$

We may also assume that there is a unitary $s'(t)\in
M_k(M_{\af_1})$ such that
\beq\label{path1+2}
s'(0)=\af_1(h(u'')), s'(1)=u'', [h(u'')]=[h(u)]\,\,\,{\rm in}\,\,\, K_1(D)\\
\andeqn \theta_2([u])=[s]\,\,\,{\rm in}\,\,\, K_1(M_{\af_1})
\eneq
for some unitary $u''\in M_k(A).$ Define $s(t)=s'(t)\oplus
\psi_2'\circ \phi_1(u'')$ for $t\in [0,1].$ As in the proof of
\ref{l1}, we may assume that $u'=u''=u.$ Now define $\theta:
K_1(A)\to K_1(M_{\af\circ \phi_1, \phi_2})$ as follows:
$\theta([u])=[v],$ where
 \beq\label{path3}
v(t)=\begin{cases} s(2t)\,\,\,\text{if}\,\,\, t\in
[0,1/2)\\
w(2(t-1/2))\,\,\,\text{if}\,\,\, t\in [1/2,1],
\end{cases}
\eneq
Thus $\theta$ gives a \hm\, from $K_1(A)$ to $K_1(M_{\af\circ
\phi_1,\phi_2})$ such that $(\pi_0)_{*1}\circ \theta={\rm
id}_{K_1(A)}.$ Since $[\af\circ \phi_1]=[\phi_1]$ in $KK(A,A)),$
we also have
$$
[\af\circ \phi_1]=[\phi_2]\,\,\,{\rm in}\,\,\,KK(A, B\otimes Q).
$$

 Then
\beq
\hspace{-0.4in}R_{\af\circ\phi_1, \phi_2}([u])(\tau)&=&
{1\over{2\pi}}\int_0^1
\tau({dv(t)\over{dt}}v(t)^*)dt\\
&=&{1\over{2\pi}}\int_0^{1/2} \tau({ds(2t)\over{dt}}s(2t)^*dt+\\
&&{1\over{2\pi}}\int_{1/2}^{1}
\tau({dw(2(t-1/2))\over{dt}}w(2(t-1/2))^*dt\\
&=& R_{\af_1\circ \Phi\circ\phi_1,\Phi\circ \phi_1}\circ
\theta_2([u])((\phi_1)_{\sharp}(\tau))+R_{\Phi\circ\phi_1,
\phi_2}\circ \theta_1([u])(\tau)=0
\eneq
for all $\tau\in T(B\otimes Q).$

Thus ${\tilde{\eta}_{\af\circ \phi_1, \phi_2}}=0.$ It follows from
Theorem \ref{T0} that $\af\circ \phi_1$ and $\phi_2$ are
asymptotically unitarily equivalent. Again, since $K_1(B\otimes
Q)$ is divisible, as in the proof of \ref{l1}, $\af\circ \phi_1$
and $\phi_2$ are strongly asymptotically unitarily equivalent.

\end{proof}

\begin{thm}\label{T1}
Let $A$ and $B$ be two unital separable amenable simple \CA s
satisfying the UCT. Let $\p$ and $\q$ be supernatural numbers of
infinite type such that $M_\p\otimes M_\q\cong Q.$ Suppose that
$A\otimes M_{\mathfrak{p}}, A\otimes M_{\mathfrak{q}},$ $B\otimes
M_{\mathfrak{p}}$ and $B\otimes M_{\mathfrak{q}}$ have tracial
rank zero.

Let $\sigma_\p: A\otimes M_{\mathfrak{p}}\to B\otimes
M_{\mathfrak{p}}$ and $\rho_{\mathfrak{q}}: A\otimes
M_{\mathfrak{q}}\to B\otimes M_{\mathfrak{q}}$ be two unital
isomorphisms. Suppose
$$
[\sigma]=[\rho]\,\,\,  in\,\,\, KK(A\otimes Q, B\otimes Q),
$$
where $\sigma=\sigma_\p\otimes {\rm id}_{M_\q}$ and
$\rho=\rho_\q\otimes {\rm id}_{M_\p}.$

Then there is an automorphism $\af\in Aut(\sigma_\p(A\otimes
M_\p))$ such that
$$
[\af\circ \sigma_\p]=[\sigma_\p]\,\,\, in\,\,\, KL(A\otimes M_\p,
B\otimes M_\p)
$$
and $\af\circ \sigma_\p\otimes {\rm id}_{M_\q}$ is strongly
asymptotically unitarily equivalent to $\rho,$ if one of the
following holds:

{\rm (i)}\, $K_i(A\otimes M_{\mathfrak{p}})$ is torsion free
($i=0,1$),

{\rm (ii)}\, $K_i(A\otimes M_\p)=Tor(K_i(A\otimes M_\p))\oplus
G_i,$ where $G_i$ is torsion free $i=0,1.$

\end{thm}

\begin{proof}
It follows from \ref{l1} that there exists $\bt\in Aut(B\otimes
Q)$ such that $ \bt\circ \sigma $ is strongly asymptotically
unitarily equivalent to $\rho.$ Moreover, $[\bt]=[{\rm
id}_{B\otimes Q}]$ in $KK(B\otimes Q, B\otimes Q).$ Now consider
two \hm s $\sigma_{\mathfrak{p}}$ and $\bt\circ
\sigma_{\mathfrak{p}}.$ One has
$$
[\bt\circ\sigma_{\mathfrak{p}}]=[\sigma_{\mathfrak{p}}]\,\,\,{\rm
in}\,\,\,KK(A\otimes M_{\mathfrak{p}}, B\otimes Q).
$$
Since $\sigma_{\mathfrak{p}}$ is an isomorphism, it is easy to see
that $\sigma_{\sharp}: T(B\otimes Q)\to T(A\otimes
M_{\mathfrak{p}})$ is an affine homeomorphism.

In case (i),  since $K_i(A\otimes M_{\mathfrak{p}})$ is torsion
free, by applying \ref{l1} again, one obtains $\af\in
Aut(\sigma_{\mathfrak{p}}(A\otimes M_{\mathfrak{p}}))$ such that
$\af\circ\sigma_{\mathfrak{p}}$ is strongly asymptotically
unitarily equivalent to $\bt\circ \sigma_{\mathfrak{p}}.$ Note
that $\sigma(A\otimes M_{\mathfrak{p}})=B\otimes
M_{\mathfrak{p}}.$ Put $\sigma_{\mathfrak{p}}'=\af\circ
\sigma_{\mathfrak{p}}$ and let $\sigma'=\af\circ
\sigma_{\mathfrak{p}}\otimes {\rm id}_{M_\q}.$

Define $\bt\circ \sigma_{\mathfrak{p}}\otimes {\rm
id}_{M_{\mathfrak{q}}}: A\otimes M_{\mathfrak{p}}\otimes
M_{\mathfrak{q}}\to (B\otimes Q)\otimes M_{\mathfrak{q}}.$ Note that
$j: M_{\mathfrak{q}}\to M_{\mathfrak{q}}\otimes M_{\mathfrak{q}}$
defined by $a\to a\otimes 1$  is (strongly) asymptotically unitarily
equivalent to an isomorphism. It follows that $\sigma'$ is strongly
unitarily equivalent to $\bt\circ \sigma_{\mathfrak{p}}\otimes {\rm
id}_{\mathfrak{q}}.$ Since $\bt\circ {\rm
id}_{M_{\mathfrak{q}}}\otimes 1,\, 1\otimes {\rm
id}_{M_{\mathfrak{q}}}: M_{\mathfrak{q}}\to \bt(1\otimes
M_{\mathfrak{q}})\otimes M_{\mathfrak{q}}$ are strongly
asymptotically unitarily equivalent (in $\bt(1\otimes
M_{\mathfrak{q}})\otimes M_{\mathfrak{q}}$), $\bt\circ \sigma$ and
$\bt\circ\sigma_{\mathfrak{p}}\otimes {\rm id}_{M_{\mathfrak{q}}}$
are strongly asymptotically unitarily equivalent.

It follows that $\sigma'$ is strongly asymptotically unitarily
equivalent to $\bt\circ \sigma.$  Consequently $\sigma'$ is
strongly asymptotically unitarily to $\rho.$

The proof of part (ii) is exactly the same but we will apply
\ref{2L1} instead.

\end{proof}

\begin{thm}{\rm (8.1 of \cite{W})}\label{T2}
Let $A$ and $B$ be two unital separable amenable simple \CA s
satisfying the UCT.  Suppose that $A\otimes C$ and $B\otimes C$
are both of tracial rank zero for any UHF-algebras. Suppose also
that
$$
(K_0(A), K_0(A)_+, [1_A], K_1(A))\cong (K_0(B), K_0(B)_+, [1_B],
K_1(B)).
$$
 Then $A\otimes {\cal Z}\cong B\otimes {\cal Z},$ if either

 {\rm (i)} $Tor(K_0(A))$  and $Tor(K_1(A))$ miss
at least one (the same) prime order,

{\rm (ii)} or $K_i(A)=Tor(K_i(A))\oplus G_i$ for some torsion free
$G_i,$ $i=0,1.$
\end{thm}

\begin{proof}
For case (i),  suppose that $Tor(K_i(A))$ misses the prime order
$p'.$ Then it is easy to find a pair of relatively prime
supernatural integers $\p$ and $\q$ such that $K_i(A\otimes M_\p)$
is torsion free ($i=0,1$). Let
$$
\Gamma:(K_0(A), K_0(A)_+, [1_A], K_1(A))\cong (K_0(B), K_0(B)_+,
[1_B], K_1(B))
$$
be the isomorphism. Let $\kappa\in KK(A, B)$ be an element which
gives $\Gamma.$

Then there are isomorphisms $\sigma_\p: A\otimes M_\p\to B\otimes
M_\p$ and $\rho_\q: A\otimes M_\q\to B\otimes M_\q$ given by
$\kappa\otimes [{\rm id}_{M_\p}]$ and $\kappa\otimes [{\rm id}_\q].$
Since $B\otimes Q$ is divisible, it is easy to see that
$$
[\sigma_\p\otimes {\rm id}_{M_\q}]=[\rho_\q\otimes {\rm
id}_{M_\p}]\,\,\,{\rm in}\,\,\, KK(A\otimes Q, B\otimes Q).
$$
Put $\sigma=\sigma_\p\otimes {\rm id}_{M_\q}$ and
$\rho=\rho_\q\otimes {\rm id}_{M_\p}.$ Then, by \ref{T1}, $\sigma$
and $\rho$ are strongly asymptotically unitarily equivalent.
Therefore $\Gamma$ can be lifted along ${\cal Z}_{\p,\q}$ (see 4.7
of \cite{W}).

One then applies Theorem 7.1 of \cite{W}.

For case (ii), there exists a sequence of integers $\{m_k\}$ such
that
$$
M_\p=\lim_{n\to\infty}(M_{m_k}, h_k),
$$
where $h_k$ is standard amplification. Thus $(h_k)_{*i}$ is a
multiplication by $m_{k+1}/m_k.$ Note also that
$$
A\otimes M_\p=\lim_{n\to\infty}(A\otimes M_{m_k},{\rm id}_A\otimes
h_k).
$$

 It follows that
$$
K_i(A\otimes M_\p)=Tor(K_i(A\otimes M_\p))\oplus G_i',
$$
where $G_i'$ is torsion free.

One then applies part (ii) of \ref{T1} as in the proof of (i).
\end{proof}






Note that in the following statement $A$ is not assumed to be of
real rank zero,  as a  priori.

\begin{cor}\label{C2}
Let $A$ be a unital separable simple ${\cal Z}$-absorbing
ASH-algebra $A$ whose projections separate the traces. Suppose
that $K_0(A)$ has the Riesz interpolation property and
$K_0(A)/Tor(K_0(A))\not\cong \Z.$

Then $A$ has tracial rank zero and $A$ is (isomorphic to) an
AH-algebra with no dimension growth and with real rank zero if
either {\rm (i)} or {\rm (ii)} of \ref{T2} hold.
\end{cor}

\begin{proof}

For any UHF-algebra $C,$ $A\otimes C$ is  approximately divisible
and its projections separate traces. It follows from \cite{BKR} that
$A\otimes C$ has real rank zero, stable rank one and weakly
unperforated $K_0(A\otimes C).$ Note that $A\otimes C$ is also an
ASH-algebra. It follows from  a result of Winter (\cite{W2}) that
$A\otimes C$ has tracial rank zero. Since $A$ is ${\cal Z}$-absorbing,
$K_0(A)$ is weakly unperforated  (Prop. 1.2 of \cite{TW}).  Now by \cite{EG}, there exists
a unital simple AH-algebra $B$ with real rank zero and with no
dimension growth such that
$$
(K_0(A), K_0(A)_+, [1_A], K_1(A))\cong (K_0(B), K_0(B)_+, [1_B],
K_1(B)).
$$
It follows from \cite{W2} that $B\otimes {\cal Z}$ has tracial
rank zero. By \ref{T2}, $A\cong B\otimes {\cal Z}.$

\end{proof}

\vspace{0.2in}

Please see section 8 of \cite{W} for other consequences of
\ref{T2} and further discussion.






\end{document}